\documentclass{amsart}
\usepackage{amsmath,amssymb,amscd,amsthm}
\usepackage[matrix,arrow]{xy}

\newcommand{\Z}         {\mathbb Z}

\newcommand{\R}         {\mathbb R}
\newcommand{\C}         {\mathbb C}
\newcommand{\Spinc}     {{\rm Spin_{\C}}}
\newcommand{\ra}        {\rightarrow}

\newcommand{\tensor}    {\otimes}
\newcommand{\dsum}      {\oplus}
\newcommand{\Dsum}      {\bigoplus}

\newcommand{\iso}       {\cong}

\DeclareMathOperator{\coker}{coker}
\DeclareMathOperator{\aind}{a-Ind}
\DeclareMathOperator{\tind}{t-Ind}
\DeclareMathOperator{\ev}{ev}
\DeclareMathOperator{\rank}{rank}
\DeclareMathOperator{\Spin}{Spin}

\newcommand{\dbar}      {\bar{\partial}}
\newcommand{\Nbar}      {\overline{N}}
\newcommand{\dirac}     {\partial\mskip-10mu\big/}
\newcommand{\g}         {\mathfrak{g}}
\newcommand{\gs}        {\mathfrak{g}^*}

\newtheorem{thm}{Theorem}

\newtheorem{prop}{Proposition}

\title{A $K$-Theoretic Note on Geometric Quantization}
\author{David S. Metzler}
\address{Department of Mathematics, MIT, Cambridge, MA 02139}
\email{metzler@math.mit.edu}

\begin{document}

\begin{abstract}

We show that the results of the paper  
{\em Symplectic Reduction and Riemann-Roch for Circle Actions} \cite{DGMW} of
Duistermaat, Guillemin, Meinrenken and Wu can be expressed entirely 
in $K$-theory. We show that their quantization
is simply a pushforward in $K$-theory, and use Lerman's symplectic cutting 
and the localization theorem in equivariant $K$-theory to prove that
quantization commutes with reduction. Only the case where the action is
free on the zero level set of the moment map is addressed. 

\end{abstract}

\maketitle

\section{Introduction}

  In their paper, 
{\em Symplectic Reduction and Riemann-Roch for Circle Actions} \cite{DGMW},
H. Duistermaat, V. Guillemin, E. Meinrenken and S. Wu use E. Lerman's
symplectic cutting technique \cite{L} to prove that quantization commutes
with reduction, in the case of a circle action (see \cite{GS1}). 
They define the quantization of a 
compact symplectic manifold $M$ via index theory, as the index of
the $\Spinc$ operator (the ``Riemann-Roch number'') 
associated to an almost-complex structure
compatible with the symplectic structure. If $M$ is a Hamiltonian
$G$-space, this is a virtual $G$-module, i.e. an element of $R(G)$. 
In section \ref{sec:QK} we show
that this can be expressed as the pushforward $p_! L$ of the prequantum line
bundle $L$ by the map $p: M \ra *$, provided we use 
the correct $K$-orientation. This is a straightforward use of the 
index theorem; however, since the arguments in section \ref{sec:QR}
depend heavily on signs, we go through the construction slowly 
to make sure that the orientation is correct.

The rest of the paper is devoted to showing that the proof in \cite{DGMW}
can be directly translated into $K$-theory. In particular, the 
index theorem is not necessary at this point, and the crucial ingredient is
the localization theorem in equivariant $K$-theory, 
which we review in section \ref{sec:loc}. This relates the pushforward on 
the whole manifold $M$ (i.e. $Q(M)$) to the pushforwards from the fixed point
sets of the circle action. Since the fixed point sets are trivial $G$-spaces,
their equivariant $K$-theory splits up as
        $$K_G(F) \iso K(F) \tensor R(G) \iso K(F) \tensor \Z[z,z^{-1}]$$
so we can treat $K$-classes on the fixed point sets as Laurent polynomials
with coefficients in K(F). It turns out that we only need to know a few
basic facts about the $z$-dependence of these polynomials.
In section \ref{sec:cut} we 
state without proof the properties of symplectic cutting which we need;
in particular, that the reduced space $M_G$ embeds into each of the
``cut'' spaces as a component of the fixed point set. 

In section \ref{sec:QR} we use these tools to prove that quantization 
commutes with reduction (Theorem~\ref{thm:main}). 
We show that $Q(M)^G = Q(M_G)$ by equating
both to $Q(M_+)$, where $M_+$ is one of the ``cuts'' of $M$ (Props. 
\ref{prop1} and \ref{prop2}). We prove both of these propositions by
embedding the rings $K(F)[z,z^{-1}]$ into two different rings of formal
Laurent series:
        $$K(F)[[z]]_z \text{~~and~~} K(F)[[z^{-1}]]_{z^{-1}}.$$
These embeddings correspond to the limit arguments in section 2 of 
\cite{DGMW}.   

In a forthcoming paper \cite{Met} we will 
extend the ideas in section \ref{sec:QR} to generalized 
$S^{1}$ equivariant cohomology theories. That work
will generalize the result of Kalkman \cite{Kal}
on localization for manifolds with boundary.


\vspace{18pt}

To fix notation and conventions, let $G$ be a compact Lie group, 
and let $M$ be a Hamiltonian 
$G$-space with symplectic form $\omega$ and moment map 
$\phi: M \ra \gs$ (say $(M,\omega,\phi)$ for short). (In section \ref{sec:loc}
we require $G$ to be topologically cyclic; in sections
\ref{sec:cut} and \ref{sec:QR} we specialize to the case $G = S^1$.)
Choosing a Riemannian
metric $g$ gives an almost-complex structure $J$, unique up to
isotopy, by the requirement $g(v,w) = \omega(Jv,w)$.
We assume $M$ is prequantizable, with prequantum line
bundle $L$ and connection $\nabla$, and that the action of $G$ 
extends to an equivariant action on (M,L). Then the infinitesimal action
of $G$ on sections of $L$ is given by the formula of Kostant \cite{Kos}:
\begin{equation}
  \label{kostant}
  D_v(s) = \nabla_v s - i \langle \phi, v \rangle s
\end{equation}
where $v \in \g$ and $s \in \Gamma(L)$. 
Let $F$ be the fixed point set of the $G$-action, with 
connected components $F_r$: $F = \coprod_r F_r$.
Then on $F_r$ the moment map has fixed value $\phi_r$, and
Kostant's formula reduces to 
\begin{equation}
  \label{kostantfixed}
  D_v(s) = -i \langle \phi_r, v \rangle s
\end{equation}
so $\phi_r$ must be a weight vector, and the action of $G$ on $L|_{F_r}$
has weight $-\phi_r$. 
We assume in sections \ref{sec:cut} and \ref{sec:QR} that the 
action is free on the level set $Z = \phi^{-1}(0)$.
Hence the symplectic reduction $M_G := Z/G$ is well-defined
as a smooth symplectic manifold, with symplectic form $\omega_G$, and in
fact $M_G$ is prequantizable, with line bundle $L_G := L|_Z /G.$

\section{Quantization as a $K$-theoretic pushforward}
\label{sec:QK}

We recall the definition of $Q(M)$ in \cite{DGMW}. The almost 
complex structure $J$ gives a splitting of the complexified cotangent bundle
\begin{equation}
        T^* M \tensor \C = T^* M^{1,0} \dsum T^* M^{0,1}
\end{equation}
and hence a bigrading of the exterior algebra
\begin{equation}
  \label{lambi}
  \Lambda^k (T^* M \tensor \C) \iso \sum_{p+q=k} T^*M^{p,q}
\end{equation}
and of the deRham complex
\begin{equation}
        \Omega^k(M) = \sum_{p+q=k} \Omega^{p,q}(M)
\end{equation}
Define the operator $\dbar: \Omega^{0,q}(M) \ra \Omega^{0,q+1}(M)$ by
\begin{equation}
        \dbar_{0,q} = \pi^{0,q+1} \circ d_{0,q}.
\end{equation}
This gives a sequence of maps
\begin{equation}
  \label{dbarseq}
  \xymatrix{ {\cdots} \ar[r] & {\Omega^{0,q} (M)} 
        \ar[r]^{\dbar} & {\Omega^{0,q+1}(M)}
        \ar[r] & {\cdots}
  }
\end{equation}
(This is {\em not} a complex unless $M$ is a complex manifold.)

Given a Hermitian connection $\nabla$ on the prequantum line bundle $L$,
we can form an operator
\begin{eqnarray}
  \label{dbarldef}
  \dbar_L &\mskip-12mu:& \mskip-12mu \Omega^{0,q}(M,L) \ra \Omega^{0,q+1}(M,L) \\
  \dbar_L &:=& \dbar \tensor 1 + 1 \tensor (\pi^{0,1} \circ \nabla).
\end{eqnarray}
This operator has principal symbol 
\begin{equation}
  \label{dbarsym}
  \sigma(\dbar_L)(x,\alpha)(\beta) = i \alpha^{0,1} \wedge \beta
\end{equation}
where $x \in M$, $\alpha \in T^*_x M$, and 
$\beta \in T^*_x M^{0,q} \tensor L$.

We form an elliptic operator $\dirac_L$ from $\dbar_L$ by adding it to
its adjoint:
\begin{eqnarray}
  \label{diracdef}
  \dirac_L &\mskip-13mu:& \mskip-12mu \Omega^{0,\text{even}}(M,L) 
                       \ra \Omega^{0,\text{odd}}(M,L) \\
  \dirac_L &:=& \dbar_L + \dbar^*_L. 
\end{eqnarray}

The quantization is defined as 
\begin{equation}
  \label{qdef}
  Q(M) := \aind(\dirac_L) = \ker \dirac_L - \coker \dirac_L.
\end{equation}
where I have labeled this ``a-Ind'' to emphasize the analytical nature of
this definition, as opposed to the topological one I will give below.

Given an action of $G$ on $(M,L)$ we can choose $\nabla$ to be preserved
by $G$, and hence $Q(M)$ is a virtual representation of $G$.   

The Atiyah-Singer index theorem equates the analytical index with the 
topological index:
\begin{equation}
  \label{inds}
  \aind(\dirac_L) = \tind(\dirac_L).
\end{equation}
The topological index depends only on the symbol of $\dirac_L$ as an element
of $K_G(T^*M)$, which equals the symbol of $\dbar_L$ (\cite{ASi}, 522):
\begin{equation}
  \label{diracdbar}
  \tind(\dirac_L) = \tind(\dbar_L).
\end{equation}
To calculate $\tind(\dbar_L)$ we push forward the $K$-class $\sigma(\dbar_L)$
to a point, in the following manner. First, we use the Riemannian metric
$g$ on $M$ to identify $TM$ and $T^*M$. The pullback of $\sigma(\dbar_L)$ to
$TM$ is given by the complex
\begin{equation}
\label{dbarsymdiag}
  \xymatrix@R=2ex{
    {\cdots} \ar[r] & {\Lambda^q \pi_{\scriptscriptstyle TM}^* TM 
    \tensor \pi_{\scriptscriptstyle TM}^* L} \ar[r]^\sigma & 
    {\Lambda^{q+1} \pi_{\scriptscriptstyle TM}^* TM 
    \tensor \pi_{\scriptscriptstyle TM}^* L} \ar[r] & {\cdots}  \\
    & (v,u) \ar@{|->}[r] & v \wedge u}
\end{equation}
where $v \in T_x M$ and $u \in \Lambda^q T_x M \tensor L_x$. (Note that I am
considering $TM$ as a complex vector bundle, with complex structure $J$,
so the above complex $\sigma$ is an element of $K_G(TM)$.)

Now $\sigma$  is exactly the Thom isomorphism applied to the vector bundle
$L$ (\cite{ASi}, 493). We can also express the Thom isomorphism as a 
push-forward by the zero section, call it $a: M \ra TM$:
\begin{equation}
  \label{thom}
  \sigma = \mathrm{Thom}_{\scriptscriptstyle TM} (L) = a_! L \in K_G(TM).
\end{equation}
The next step in calculating the index is to embed $M$ equivariantly in a 
trivial $G$-space $\C^n$. We can in fact choose $n$ large enough so that
the normal bundle $N_M$ will have a (unique) complex structure \cite{AH},
defined by the exact sequence of complex vector bundles
\begin{equation}
\label{ses}
  \xymatrix{ 0 \ar[r] & {TM} \ar[r]^{Tj} & {j^* T\C^n} \ar[r] &
             {N_M} \ar[r] & 0}
\end{equation}

We then have the following diagram:
\begin{equation}
  \label{pushdiag}
  \xymatrix{F \ar[r]^i & M \ar[r]^j \ar[d]_a
        & {\C^n} \ar[d]_b & {*} \ar[l]_k \ar[ld]^l \\
        & TM \ar[r]^{Tj} & {T\C^n}}
\end{equation}
Here $i$ is the inclusion of the fixed point set $F$ (which will come into 
the picture soon), $j$ is the chosen embedding of $M$ in $\C^n$, with 
corresponding tangent map $Tj$; $a$ and $b$ are the zero sections of $TM$
and $T\C^n$ respectively; and $k$ and $l$ are the inclusions of the origin
into $\C^n$ and $T\C^n \iso \C^n \dsum \C^n$ respectively. 

The topological index of $\dbar_L$ is defined to be the pushforward
\begin{eqnarray}
  \label{tinddef}
  \tind(\dbar_L) &=& (l_!)^{-1} (Tj)_! \sigma \\
                 &=& (l_!)^{-1} (Tj)_! a_! L.
\end{eqnarray}
It would seem that functoriality of the pushforward immediately gives
\begin{equation}
  \label{qisp}
  Q(M) = \tind(\dbar_L) = (k_!)^{-1} j_! L =: p_! L
\end{equation}
where the last equality is the definition of the pushforward of $L$ by
the map $p: M \ra *$. However, we need to be careful about $K$-orientations.
Each of these pushforwards requires a $K$-orientation, for example a
complex structure, on the corresponding normal bundle for 
a precise definition. Let us adopt the temporary notation $\kappa(f)$
to denote the complex structure on the normal bundle to an embedding
$f: X \ra Y$ to be used in the pushforward $f_!: K_G(X) \ra K_G(Y)$.
Then in the diagram (\ref{pushdiag}) we know 
\begin{eqnarray*}
  \label{kappas}
  \kappa(a)   &=& TM \\
  \kappa(b)   &=& T\C^n = \underline{\C^n} \qquad 
                    \hbox{(trivial rank n bundle)}\\
  \kappa(Tj)  &=& \pi_{\scriptscriptstyle TM}^* N_M \dsum 
                    \pi_{\scriptscriptstyle TM}^* \Nbar_M  \\
  \kappa(k)   &=& \underline{\C^n} \\
  \kappa(l)   &=& \underline{\C^{2n}} = \kappa(k) \dsum k^* \kappa(b)
\end{eqnarray*}
(See \cite{ASi} for the identification of $\kappa(Tj)$.) The last equation
shows that the right-hand triangle in the pushforward diagram 
\begin{equation}
  \label{kpushdiag}
  \xymatrix{{K(F)} \ar[r]^{i_!} & {K_G(M)} \ar[r]^{j_!} \ar[d]_{a_!}
        & {K_G(\C^n)} \ar[d]_{b_!} \ar[r]^{(k_!)^{-1}} & {K_G(*)} \\   
        & {K_G(TM)} \ar[r]^{(Tj)_!} & {K_G(T\C^n)} \ar[ru]_{(l_!)^{-1}}}
\end{equation}
commutes (this is just Bott periodicity); what we need is the correct 
$\kappa(j)$ to make the square commute. But if the square is to commute
we must have
\begin{eqnarray*}
  \label{kappaj}
  j^* \kappa(b) \dsum \kappa(j) &=& \kappa(Tj \circ a)  \\ 
                     &=& \kappa(a) \dsum a^* \kappa(Tj) \\
                     &=& TM \dsum a^* (\pi^* N_M \dsum \pi^* \Nbar_M) \\
                     &=& TM \dsum N_M \dsum \Nbar_M \\
                     &=& \underline{\C^n} \dsum \Nbar_M \\
                     &=& j^* \kappa(b) \dsum \Nbar_M 
\end{eqnarray*}
so we must choose $\kappa(j) = \Nbar_M$, and not $N_M$. 
This gives us our
\begin{thm}
\label{thm:qisp}
  Let $p: M \ra *$ be the unique map and let $N_M$ be the stable 
  normal bundle to $M$, defined by (\ref{ses}) above. Then the quantization
  of (M,L) is exactly $p_! L$, using the $K$-orientation $\Nbar_M$.  
\end{thm}

The fact that we must use $\Nbar_M$ will be significant when we look at 
localization in section \ref{sec:loc}. There we will be concerned with the
fixed point set $F$ of the $G$-action. We have a diagram
\begin{equation}
\label{tridiag}
  \xymatrix{K_G(F) \ar[r]^{i_!} \ar[dr]_{q_!} & K_G(M) \ar[d]^{p_!} \\
            & K_G(*)  }  
\end{equation}
which we want to commute. We now know how to precisely define $p_!$ and
$q_!$ to agree with the quantization: we use the complex structures 
$\kappa(p) = \Nbar_M$ and $\kappa(q) = \Nbar_F$ on the respective
stable normal bundles. Letting $N$ be the usual complex normal bundle of 
$F$ in $M$, defined by 
\begin{equation}
\label{MFses}
  \xymatrix{ 0 \ar[r] & {TF} \ar[r]^{Ti} & {i^* TM} \ar[r] &
             N \ar[r] & 0}
\end{equation}
we must have
\begin{eqnarray}
  \label{MFkappa}
  \Nbar \dsum i^* \Nbar_M &=& \Nbar_F   \\
                              &=& \kappa(q)   \\
                              &=& \kappa(i) \dsum i^* \kappa(p) \\
                              &=& \kappa(i) \dsum i^* \Nbar_M
\end{eqnarray}
so $\kappa(i) = \Nbar$. This will be important in getting the signs
correct in the next section.

\section{Localization in Equivariant $K$-theory}
\label{sec:loc}
  The key tool we use is the localization theorem of Atiyah and Segal 
\cite{ASe} \cite{S} in equivariant $K$-theory, which we briefly review. 
We follow the treatment in \cite{ASe} except that they are doing index theory 
and hence work on the tangent bundle, while we work on $M$ itself. 

We wish to calculate $Q(M) = p_! L \in K_G(*) = R(G)$. Since every element
of $R(G)$ is determined by its character, we can specify $Q(M)$ by evaluating
its character at every element $g \in G$, or even on a dense subset of 
elements $g$. For simplicity, assume $G$ is topologically cyclic. (Of course
eventually $G$ will simply be $S^1$.) Fix a (topological)  generator
$g \in G$, i.e. let $(g)$ be dense in $G$. Then $M^g = M^G = F$. 
The localization theorem gives a formula for computing the character of
$p_! L$, evaluated at $g$, in terms of data on $F$. 

We start with the diagram (\ref{tridiag}). All of these rings are actually
$R(G)$-algebras, so we can localize at $g$ (this inverts all characters
not vanishing at $g$).
\begin{equation}
  \label{triloc}
  \xymatrix{K_G(F)_g \ar[r]^{(i_!)_g} \ar[dr]_{(q_!)_g} 
            & K_G(M)_g \ar[d]^{(p_!)_g} \\
            & R(G)_g  }  
\end{equation}
\begin{thm}[\protect\cite{ASe}]
\label{thm:loc}
    The map $(i_!)_g$ is an isomorphism of $R(G)_g$-modules. 
\end{thm}
This allows us to calculate the pushforward by $p$ in terms of the
pushforward by $q$, at least in the localized ring $R(G)_g$. This is good
enough, since we are interested in evaluating $p_! L$ at $g$, and the
evaluation map $\ev_g: R(G) \ra \C$ factors through the localization 
$R(G)_g$. In fact we have the following commutative diagram:
\begin{equation}
  \label{bigdiag}
  \xymatrix{
  {K_G(M)} \ar[d]_{p_!} \ar[r] 
  & {K_G(M)_g} \ar[d]_{(p_!)_g} \ar[r]^{(i_!)_g^{-1}}_{\sim}
  & {K_G(F)_g}  \ar[ld]^{(q_!)_g} \ar[r]^-{\sim}
  & {K(F) \tensor R(G)_g} \ar[d]_{1 \tensor \ev_g}\\
  {R(G)} \ar[r] \ar[rd]_{\ev_g} 
  & {R(G)_g} \ar[d]_{\ev_g} & & {K(F) \tensor \C} \ar[dll]^{q_!}\\
  & {\C}    }  
\end{equation}
Here we have used the isomorphism $K_G(F) \iso K(F) \tensor R(G)$ for
the trivial $G$-space $F$. 

The next step is to explicitly identify the map $(i_!)_g^{-1}$. 
This turns out to be simple. Recall that the $K$-orientation for
the map $i$ was $\kappa(i) = \Nbar$. We have (\cite{ASi}, 493)
\begin{equation}
  \label{starshriek}
  i^* i_! u = \Lambda \Nbar \cdot u := \left( \sum (-1)^k 
               \Lambda^k \Nbar \right) \cdot u  
\end{equation}
so when we localize at $g$, the inverse is simply
\begin{equation}
  \label{inv}
  (i_!)_g^{-1} L = \frac{i^* L}{\Lambda \Nbar}.
\end{equation}
Using this explicit inverse we can write down the localization formula
giving the result of evaluation at $G$:
\begin{equation}
  \label{locform}
  (p_! L)(g) = q_! \left( \frac{i^* L (g)}{\Lambda \Nbar (g)} \right)
\end{equation}
where the quantity in parentheses is in $K(F) \tensor \C$, 
and the evaluations
$i^*L(g)$ and $\Lambda \Nbar (g)$ are defined by the composite map 
$$\xymatrix{K_G(F)_g \iso K(F) \tensor R(G)_g \ar[r]^-{1 \tensor \ev_g}
                & K(F) \tensor \C }.$$
In other words, to use this formula, we need to represent $i^* L$ and
$\Lambda \Nbar (g)$ as sums of $G$-fixed bundles tensored with characters of
$G$, and then evaluate at $g$ by the prescriptions
\begin{eqnarray}
  \label{eval}
  u \tensor \chi        &\mapsto& u \cdot \chi(g)             \\
  u \tensor \chi/\psi   &\mapsto& u \cdot \chi(g)/\psi(g).
\end{eqnarray}

\section{Quantization Commutes With Reduction}
\label{sec:QR}

\subsection{Symplectic Cutting}
\label{sec:cut}

From here on we deal only with the case $G = S^1$.
In \cite{L} E. Lerman defines an operation 
on a Hamiltonian ${S^1}$-space
called {\em symplectic cutting}. (See also \cite{DGMW}.) 
Cutting produces from a Hamiltonian ${S^1}$-space $(M,\omega,\phi)$ 
a pair of Hamiltonian ${S^1}$-spaces $(M_+,\omega_+,\phi_+)$ 
and $(M_-,\omega_-,\phi_-)$ with the following properties:
\begin{itemize}
  \item The reduced space $M_{{S^1}}$ embeds in both $M_+$ and $M_-$ 
        as a component of the fixed point set.
  \item $M_+ \setminus M_{{S^1}}$ is equivariantly, symplectically isomorphic
        to $\phi^{-1}(\R_+)$.
  \item $M_- \setminus M_{{S^1}}$ is equivariantly, symplectically isomorphic
        to $\phi^{-1}(\R_-)$.  
  \item $\phi_+(M_{{S^1}}) = \phi_-(M_{{S^1}}) = 0$. 
\end{itemize}
Further, if $M$ is prequantizable, with prequantum line bundle $L$ and
prequantizing connection $\nabla$, both $M_+$ and $M_-$ are prequantizable,
with line bundles $L_+$ and $L_-$, and the restriction of these bundles
to the reduced space is just the reduced line bundle:
\begin{equation}
  \label{lbcompat}
  L_+|_{M_{{S^1}}} \iso L_{{S^1}}, \qquad L_-|_{M_{{S^1}}} \iso L_{{S^1}}.
\end{equation}


\subsection{The Main Results}
\label{subsec:setup}

Since symplectic cutting embeds the reduced space $M_{S^1}$ into pieces of
the original space $M$ as a fixed point set, we can apply the 
$K_{S^1}$-localization theorem. We can prove our main result, 
Theorem~\ref{thm:main} by comparing $M$, $M_+$, and $M_{S^1}$.

In sections \ref{subsec:laurent} and \ref{subsec:proof} 
we use Laurent series expansions to prove the
following two propositions. 

\begin{prop}
  \label{prop1}
  Let $M$, $N$ be prequantizable Hamiltonian $S^1$-spaces with 
  moment maps $\phi$, $\psi$. 
  Assume that $0$ is not the maximum value of $\phi$ or of $\psi$.
  If $\phi^{-1}(\R_+)$ is equivariantly symplectomorphic to $\psi^{-1}(\R_+)$,
  then
        $$Q(M)^{S^1} = Q(N)^{S^1.}$$
\end{prop}

\begin{prop}
  \label{prop2}
  Let $M$ be a prequantizable Hamiltonian $S^1$-space with moment map $\phi$
  and line bundle $L$.
  Assume that $0$ is the minimum value of $\phi$. 
  Let $F_0 = \phi^{-1}(0)$ and consider the maps 
  $i_0: F_0 \ra M$, $q_0: F_0 \ra *$. Then
        $$Q(M)^{S^1} = (q_0)_! i_0^* L \in K(*) \iso \Z.$$
\end{prop}

Assuming these propositions we can prove that quantization commutes with
reduction. Consider our Hamiltonian $S^1$-space M.
Since $0$ is a regular value of $\phi$ it is certainly not the maximum of
$\phi$ or of $\phi_+$. Applying Prop.~\ref{prop1} to $M$ and $M_+$ gives
\begin{equation}
  \label{MM+}
  Q(M)^{S^1} = Q(M_+)^{S^1}.
\end{equation}

Now $0$ {\em is} the minimum value of $\phi_+$, so we can apply 
Prop.~\ref{prop2} to $M$ and $\phi^{-1}(0) = M_{S^1}$. Here
$i_0^* L = L_{S^1}$ and $(q_0)_! L_{S^1} = Q(M_{S^1})$, so
\begin{equation}
  \label{M+MG}
  Q(M_+)^{S^1} = Q(M_{S^1}).
\end{equation}

Putting these together gives our main theorem.

\begin{thm}
\label{thm:main}
  Let $(M,\omega,\phi)$ be a prequantizable Hamiltonian $S^1$-space with
  prequantum line bundle $L$ and assume that
  the action is free on the zero level set $\phi^{-1}(0)$. 
  Then the quantization $Q(M) = p_!(L)$ commutes with reduction:
        $$Q(M)^{S^1} = Q(M_{S^1}).$$
\end{thm}

\vspace{8mm}

\subsection{Expansion in Laurent Series}
\label{subsec:laurent}


It remains to prove Props. \ref{prop1} and \ref{prop2}. We will use the
localization theorem, and two different expansions in Laurent Series,
which correspond to the limit arguments (``$z \ra 0$'' and ``$z \ra \infty$'')
in \cite{DGMW}.

In the case of a Hamiltonian circle action, we can express the localization
formula~(\ref{locform}) in the following way. First, we recall that the
fixed point set $F$ breaks up into connected components $F_r$, on each of
which the action of ${S^1}$ on $L$ has weight $-\phi_r$. The localization
formula becomes
\begin{equation}
  \label{loc_r}
  (p_! L)(g) = q_! \left( \sum_r \frac{i_r^* L (g)}
                                        {\Lambda \Nbar_r (g)} \right)
\end{equation}

It turns out that we need to know very little about the quantity $p_! L$
to prove the propositions. This allows us to do everything within 
$K_G(F)_g \tensor \C \iso K(F) \tensor \C[z,z^{-1}]_g$, without 
actually evaluating the pushforward. 
(We tensor with $\C$ so that
later operations involving tensors will be exact; since the final
result $(p_! L)(g)$ is in $\C$ this is sufficient.)
  
We want to consider the contribution of each component of the fixed point
set in turn. So fix $r$, and let $l_r = i_r^* L$ with the {\em trivial} action 
of $G$. Then as an element of 
        $$K_G(F_r)_g \tensor \C \iso K(F_r) \tensor R(G)_g \tensor \C 
                \iso K(F_r) \tensor \C[z,z^{-1}]_g$$
we have
\begin{equation}
  \label{l_r}
  i_r^* L = l_r \cdot z^{-\phi_r}.
\end{equation}
Let $S_r$ be the set of weights of the action of $G$ on $N_r$.
Then we can write $\Nbar_r \in K(F_r) \tensor \C[z,z^{-1}]_g$ as
\begin{equation}
  \label{N_rk}
  \Nbar_r = \Dsum_{k \in S_r} \Nbar_{r,k} z^{-k}
\end{equation}
where the $\Nbar_{r,k}$ are vector bundles with trivial $G$-action. (Note they
are not necessarily line bundles. In fact we will not need to use a
splitting principle.) Let $n_{r,k} = \rank \Nbar_{r,k}$. 

The only difficult step is inverting $\Nbar_r$. To formally invert
polynomials, it is useful to embed the polynomial ring in the
larger ring of formal power series. In our case, we need to embed 
our ring 
        $$K(F) \tensor \C[z,z^{-1}]_g$$
of formal Laurent polynomials with coefficients in $K(F)$ 
(localized at $g$), 
into two different rings of formal Laurent series, depending on which 
proposition we want to prove:
\begin{eqnarray}
  \label{twoembeds}
  K(F) \tensor \C[z,z^{-1}]_g &\ra& K(F) \tensor \C[[z]]_z        \\
  K(F) \tensor \C[z,z^{-1}]_g &\ra& K(F) \tensor \C[[z^{-1}]]_{z^{-1}}.
\end{eqnarray}
Here $\C[[z]]_{z}$ is the ring of Laurent series ``at $z=0$,'' i.e.
allowing an infinite number of nonzero terms with positive powers
of $z$, but only a finite number of nonzero terms with negative powers
of $z$. Similarly $\C[[z^{-1}]]_{z^{-1}}$ is the ring
of Laurent series ``at $z=\infty$.''
It is not hard to see that these maps really are injective, since
they are derived from the basic inclusion $\C[z] \subset \C[[z]]$ 
by localization and tensoring, and both operations are exact over a field.
Hence we lose no information in this process.

Let $S_{r+} := S_r \cap \Z_+$ and $S_{r-} := S_r \cap \Z_-$. 
Note $S_r = S_{r+} \cup S_{r-}$ since the zero weight doesn't appear in the
normal bundle. 
Also note that $S_{r+} = \emptyset$ iff $\phi_r$ is the maximum of $\phi$,
and $S_{r-} = \emptyset$ iff $\phi_r$ is the minimum of $\phi$. 
(For a general manifold, these would only be statements about local
minima and maxima, but since this is a Hamiltonian $S^1$-space there are
no local maxima or minima except the global max and min. 
See \cite{GS:STP}, \cite{Kir}.)
Then
\begin{eqnarray*}
  \Lambda \Nbar_r &=& \Lambda \Dsum_{k \in S_r} \Nbar_{r,k} z^{-k}         \\
                &=& \prod_{k \in S_r} \Lambda (\Nbar_{r,k} z^{-k})       \\
                &=& \prod_{k \in S_r} \sum_{j=0}^{n_{r,k}} (-1)^j 
                        (\Lambda^j \Nbar_{r,k}) z^{-jk}              \\
                &=& P_r(z) Q_r(z^{-1})
\end{eqnarray*}
where $P_r$ and $Q_r$ are polynomials with constant term $1$ and invertible
leading coefficient:
\begin{eqnarray*}
  \text{Leading coeff. of~} P_r &=& \prod_{k \in S_{r-}} 
                                (-1)^{n_{r,k}} \det \Nbar_{r,k}         \\
  \text{Leading coeff. of~} Q_r &=& \prod_{k \in S_{r+}} 
                                (-1)^{n_{r,k}} \det \Nbar_{r,k}.         
\end{eqnarray*}
We have 
\begin{eqnarray*}
  P_r = 1 \!\!&\iff&\!\! S_{r-} = \emptyset \iff \phi_r = \phi_{min}     \\
  Q_r = 1 \!\!&\iff&\!\! S_{r+} = \emptyset \iff \phi_r = \phi_{max}.
\end{eqnarray*}
Hence we can invert $\Lambda \Nbar$ in the formal Laurent
rings according to the results in the appendix. 

\begin{enumerate} 
  \item (``Limit as $z \ra 0$'') In the ring $K(F) \tensor \C[[z]]_z$
    we have 
\begin{eqnarray}
  (\Lambda \Nbar_r)^{-1} &=& P_r^{-1} Q_r^{-1}                        \\
\label{pospowers}        &=& \left\{ 
                             \begin{array}{ll}
                                1 + O(z) & \text{~if~} 
                                                \phi_r = \phi_{max}   \\
                                O(z)       & \text{~if~}
                                                \phi_r \ne \phi_{max}. 
                             \end{array} \right.
\end{eqnarray}
Here $O(z^k)$ indicates a term that has no nonzero coefficients below
the $k$th power.

  \item (``Limit as $z \ra \infty$'') In $K(F) \tensor \C[[z^{-1}]]_{z^{-1}}$
    we have 
\begin{eqnarray}
  (\Lambda \Nbar)^{-1} &=& P^{-1} Q^{-1}                        \\
  \label{negpowers}    &=& \left\{ 
                             \begin{array}{ll}
                                1 + o(z^{-1}) & \text{~if~} 
                                                   \phi_r = \phi_{min}  \\
                                o(z^{-1})       & \text{~if~} 
                                                   \phi_r \ne \phi_{min}.
                             \end{array} \right.
\end{eqnarray}
Here $o(z^k)$ indicates a term that has no nonzero coefficients {\em above}
the $k$th power.
\end{enumerate}

\subsection{Proof of the Propositions}
\label{subsec:proof}

\begin{proof}[Proof of Prop. 1.] 

The multiplicity of the trivial representation in $Q(M)$ 
is just the constant term in the polynomial 
\begin{equation}
  \label{QMz}
  Q(M)(z) =  q_! \left( \sum_r l_r z^{-\phi_r} 
                       P_r^{-1} Q_r^{-1} \right).
\end{equation}
We will show that the constant term depends only on the 
components $F_r$ with $\phi_r > 0$.

We can express 
\begin{equation}
  \label{QMK}
  \sum_r l_r z^{-\phi_r} P_r^{-1} Q_r^{-1}
\end{equation}
in the Laurent series ring $K(F){[[z]]_z}$ using (\ref{pospowers}). 
The terms in the
sum~(\ref{QMK}) with $\phi_r < 0$ are of the form $O(z)$ by (\ref{pospowers}), 
so they do not contribute to the constant term. 
The terms with $\phi_r = 0$ are also of the form $O(z)$ since we are assuming
$\phi_{max} \ne 0$. 

Since the constant terms in $Q(M)(z)$ and $Q(N)(z)$ only depend on the
fixed point sets in $\phi^{-1}(\R_+)$ and $\psi^{-1}(\R_+)$ respectively, 
and these portions of $M$ and $N$ are symplectomorphic, we have 
        $$Q(M)^{S^1} = Q(N)^{S^1}.$$
\end{proof}

\begin{proof}[Proof of Prop. 2.] 
$(M,\omega,\phi)$ is a prequantizable Hamiltonian $S^1$-space 
with line bundle $L$. Since $0$ is the minimum value of $\phi$, 
$F_0 = \phi^{-1}(0) \subset M^{S^1}$. The localization theorem gives
\begin{equation}
  \label{QM+z}
  Q(M)(z) =  q_! \left( \sum_{\phi_r > 0} l_r z^{-\phi_r} 
                P_r^{-1} Q_r^{-1} \right)
                + (q_0)_! (i_0^* L P_0^{-1} Q_0^{-1}) 
\end{equation}
where $q_0: F \ra *$ and $i_0: F \ra M$. 
 
In $K(F){[[z^{-1}]]_{z^{-1}}}$, the terms in the summation 
\begin{equation}
  \label{QM+K}
  \sum_{\phi_r > 0} l_r z^{-\phi_r} P_r^{-1} Q_r^{-1}
\end{equation}
have only negative powers, by (\ref{negpowers}), so they do not contribute to
the constant term. The contribution of F is 
        $$i_0^* L \cdot (1 + o(z^{-1}))$$
again by (\ref{negpowers}), so the constant term is just the pushforward
from $F_0$,
\begin{equation}
  \label{last}
  Q(M)^{S^1} = (q_0)_! i_0^* L.
\end{equation}
\end{proof}

\section{Appendix: Inverting Polynomials in Laurent Series Rings}
\label{sec:laurapp}

Here we write down some elementary facts about inverting polynomials
in rings of formal Laurent Series, needed in section \ref{subsec:laurent}
above. 

Let $R$ be a ring. We want to formally invert Laurent polynomials, i.e.
elements of $R[z,z^{-1}]$. For our purposes we only need to know the most
basic facts about the dependence of these inverses on $z$, and for that
purpose, it is useful to embed $R[z,z^{-1}]$ into the two rings of 
formal Laurent series, $R[[z]]_z$ and $R[[z^{-1}]]_{z^{-1}}$. 

We look at the case of $R[[z]]_z$, formal Laurent
series in at $z=0$. 

\begin{enumerate}
  \item Let $a(z) = a_0 + a_1 z + \ldots + a_n z^n$ be a polynomial over 
        $R$. Suppose $a_0$ is invertible. Then we can invert $a$ in
        $R[[z]]$, hence {\em a fortiori} in $R[[z]]_z$:
        \begin{eqnarray*}
          a^{-1} &=& a_0^{-1} (1 + (a_1/a_0) z + \ldots + (a_n/a_0) z^n)^{-1}\\
                 &=& a_0^{-1} \sum_{l=0}^\infty 
                        ((a_1/a_0) z + \ldots + (a_n/a_0) z^n)^l \\
                 &=& a_0^{-1} + O(z).
        \end{eqnarray*}

  \item Let $b(z) = b_0 + b_1 z^{-1} + \ldots + b_m z^{-m}$ 
        be a Laurent polynomial over $R$. Suppose $b_m$ is invertible. 
        Then we can invert $b$ in $R[[z]]_z$:
        \begin{eqnarray*}
          b      &=& b_m z^{-m} (1 + (b_{m-1}/b_m) z + \ldots + 
                        (b_0/b_m) z^m) \\
          b^{-1} &=& b_m^{-1} z^m \sum_{l=0}^\infty 
                        ((b_{m-1}/b_m) z + \ldots + (b_0/b_m) z^m)^l \\
                 &=& b_0^{-1} z^m + O(z^{m+1}).
        \end{eqnarray*}
\end{enumerate}

The case of $R[[z^{-1}]]_{z^{-1}}$ is exactly parallel; simply
exchange $z$ with $z^{-1}$ and $O$ with $o$.

\end{document}